\documentclass[oneside,a4paper]{article}
\usepackage{latexsym}
\usepackage{graphicx}
\usepackage{epsfig}
\usepackage[cp1250]{inputenc}
\begin{document}

\setcounter{page}{15}

\begin{center}
\vspace*{30mm}

\normalsize {\textbf{STATE-SPACE CONTROLLER DESIGN FOR THE
                     FRACTIONAL-ORDER REGULATED SYSTEM}}   \\

\vspace{5mm}
\normalsize
\textbf{\v{L}ubom\'{\i}r DOR\v{C}\'AK$^{\mbox{1}}$},
\textbf{Ivo PETR\'A\v{S}$^{\mbox{1}}$}  \\
\textbf{Imrich KO\v{S}TIAL$^{\mbox{1}}$},
\textbf{J\'an TERP\'AK$^{\mbox{1}}$}  \\
$^{\mbox{1}}$
Department of Informatics and Process Control \\
BERG Faculty, Technical University of Ko\v{s}ice \\
B. N\v{e}mcovej 3, 042 00 Ko\v{s}ice, Slovak Republic \\
e-mail:
\textit{\{lubomir.dorcak,ivo.petras,imrich.kostial,jan.terpak\}@tuke.sk} \\
\end{center}
\vspace{3mm} 

\hspace{-6.5mm}
\textbf{Abstract:}
 In this paper we will present a mathematical description and analysis
of a fractional-order regulated system in the state space and the
state-space controller design based on placing the closed-loop poles on
the complex plane. Presented are the results of simulations
and stability investigation
of this system.
\vspace{3.0mm}

\hspace{-6.5mm}
\vspace{-2.0mm}

\section*{\normalsize{1. \hspace{1.2mm}Introduction}}
\hspace{5.0mm}
Research of the \, fractional-order \, derivatives continued in the last
decades not only as a \, mathematical branch \, \cite{Oldham, Samko,
Podlubny6, PodlubnyMat}, \, etc., \, but in many applied domains, \, e.g.
\cite{Westerlund, Podlubny5} and in control theory too
\cite{Outstaloup} -
\cite{Dorcak}, etc.
In works \cite{Outstaloup, Axtell, Kalojanov} the first generalizations
of analysis methods for fractional-order control systems were made
($s$-plane, frequency response, etc.). The following of the
above-mentioned works were oriented to the methods of fractional-order
system parameters identification, methods of fractional-order
controllers synthesis, methods of stability analysis, methods of control
of chaotic fractional-order systems, and so on.

 In works \cite{Matignon, Vinagre1} was presented a state space
model described in vector and matrix relations expressing the
fractional-order derivatives
\begin{eqnarray} \label{r1}
        {\bf x}^{(\alpha)}(t) =
           {\bf A}~{\bf x}(t) + {\bf B}~u(t), \nonumber \\
         y(t) = {\bf C}~{\bf x}(t),
   \hspace{1.5em} t \geq 0. \hspace{0.2em}
\end{eqnarray}
This description is convenient only for simple models of systems
with only one fractional-order derivation, or for spatial type of
systems.
In work \cite{Dorcak} we proposed a state space model of the linear
time-invariant one dimensional system which expresses the first
derivatives in the state space equations and has the \,
{\it {classical \, state \, space interpretation}}
for the fractional-order system too. On the
right side of these equations we can then transfer more than one
fractional-order derivatives of the state space variables. A disadvantage
of this expression is that in time domain we cannot express the state
space equations in vector and matrix relations as in
 previous
description. But we can do this in $s$-plane.

This contribution deals with a mathematical description and analysis
of a fractional-order regulated system in the state space and the
state-space controller design based on placing the closed-loop poles
on the complex plane.

\vspace{-1.0mm}

\section*{\normalsize{2. \hspace{1.2mm}Definition of the system}}
\vspace{-2.0mm}

\hspace{5.0mm}
        For the definition of the control system we consider
a simple unity feedback control system illustrated in Fig.1,
where $G_{s}(s)$ denotes the transfer function of the controlled
system and $G_{r}(s)$ is the controller transfer function, both
integer- or fractional-order.
\begin{figure}[ht]
\centering
\includegraphics[scale=0.68]{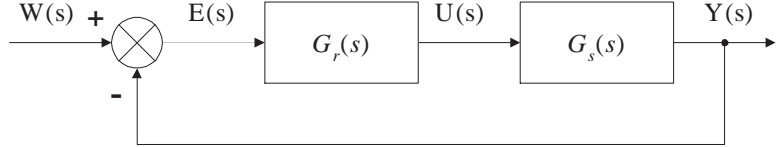}
\vspace{-2.0mm}
\caption{Unity-feedback control system}
         \label{fig:FBCS}
\end{figure}
        The differential equation of the above closed regulation
system for the transfer function of the controlled system
$G_{s}(s)=1/(a_2s^{\alpha}+a_1s^{\beta}+a_0)$
and the fractional-order $PD^{\delta}$ controller
$G_{r}(s)=K+T_ds^{\delta}$
has the form
\vspace{-0.8mm}
\begin{equation} \label{r2}
     a_{2}\, y^{(\alpha)}(t) +
     a_{1}\, y^{(\beta)}(t) +
     T_{d}\, y^{(\delta)}(t) +
     (a_{0}+K)\, y(t) =
     K\, w(t) +
     T_{d}\, w^{(\delta)}(t)
\end{equation}
\vspace{-2.2mm}
where  $\alpha, \beta, \delta$ are generally  real
numbers and $a_0, a_1, a_2, K, T_d$ are  arbitrary
constants.

\vspace{1.0mm}

\section*{\normalsize{3. \hspace{1.2mm}Fractional-order~control~system
                          with $PD^{\delta}$~controller~in~the~state~space}}
\vspace{-1.0mm}

\hspace{5.0mm}
 Consider the system described by differential equation (\ref{r2}).
After its modification and substitution of state space variables
$x(t) = x_1(t), \dot{x}(t) = x_2(t), \ddot{x}(t) = \dot{x_2}(t)$
we can derive the following state space model equivalent to model
(\ref{r2})
\vspace{-2.0mm}
\begin{eqnarray} \label{r4}
    \dot{x_1}(t) = x_2(t)~,
     \hspace{10.1cm}
     \nonumber \\
    \dot{x_2}(t) = - \frac{a_0+K}{a_2}~x_1^{(2-\alpha)}(t)
                   - \frac{T_d}  {a_2}~x_2^{(1+\delta-\alpha)}(t)
                   - \frac{a_1}  {a_2}~x_2^{(1+\beta-\alpha)}(t)
                   + \frac{1}    {a_2}~w^{(2-\alpha)}(t)~,
     \hspace{0em} \\
        y  (t) =   K~x_1(t)
                 + T_d~x_2^{(\delta-1)}(t)~,
     \hspace{1.5em} t \geq 0. \hspace{17.1em}
     \nonumber
\vspace{-1.0mm}
\end{eqnarray}
Of course, we can make other alternative state space models for
the same system.
For the approximation of the fractional-order derivatives on the
right-hand side of equations (\ref{r4}) we can take the relation
(\ref{rr4}) from e.g. \cite{Samko, Podlubny6, Dorcak2}
\vspace{-2mm}
\begin{equation} \label{rr4}
   y^{(\alpha)}(t) \approx
   h^{-\alpha} \sum_{j=0}^{N(t)}b_{j}y(t-j h ) \; ,
\end{equation}
where $L$ is "memory length", $h$ is
 time step of calculation,
$N(t) = \min \left\{
           \left[
              \frac{\mbox{$t$}}{\mbox{$h$}}
           \right] , \;
           \left[
              \frac{\mbox{$L$}}{\mbox{$h$}}
           \right]
         \right\} , $
$[z]$ is the integer part of $z$,
$b_j = (-1)^j {{\alpha} \choose j}$,
and $ {{\alpha} \choose j} $ is binomial coefficient. To calculate
$b_j$ it is convenient to use the following recurrent relation
\begin{equation} \label{rr6}
       b_0 = 1
         \ , \ \
       b_j = (1- \frac {1+\alpha} {j}) \ b_{j-1}
\end{equation}
  After the above-mentioned fractional-order derivatives discretisation
(\ref{rr4}) and discretisation of the first derivatives on the left-hand
side of equations (\ref{r4}) we obtained the simple Euler methods for
solving the state space model
\vspace{-1mm}
\begin{eqnarray} \label{rrr4}
    x_{1,k+1} = x_{1,k} + h \, x_{2,k}~,
     \hspace{8.0cm} \nonumber \\
    x_{2,k+1} = x_{2,k} + h \, \Bigl(
                -\frac{a_0+K}{a_2}~h^{\alpha-2}
                                      \sum_{j=0}^{k}b_j x_{1,k-j}
                -\frac{T_d}  {a_2}~h^{\alpha-\delta-1}
                                      \sum_{j=0}^{k}c_j x_{2,k-j}
     \hspace{0.5cm} \nonumber \\
                -\frac{a_1}  {a_2}~h^{\alpha-\beta-1} \hspace{0em}
                                      \sum_{j=0}^{k}d_j x_{2,k-j}
                 \hspace{0.2em} + \hspace{0.3em}
                 \frac{1}    {a_2}~h^{\alpha-2} \hspace{0.5em}
                                      \sum_{j=0}^{k}b_j w_{k-j}
                            \, \Bigr)
     \hspace{1.2em} \\
        y_k  =   K~x_{1,k} +
                 T_d~h^{(1-\delta)}
                                      \sum_{j=0}^{k}e_j x_{2,k-j}
     \hspace{3.5em} k \geq 0. \hspace{9.0em}
     \nonumber
\end{eqnarray}
From these equations we can compute state trajectories of the
fractional-order control system described above.

\section*{\normalsize{4. \hspace{1.2mm}Design of the $PD^{\delta}$
                                       controller~in~the~state~space}}

\hspace{5.0mm}
After Laplace transformation we can write state space equations (\ref{r4})
in vector and matrix relations in $s$-plane
\vspace{-1mm}
\begin{eqnarray} \label{r5}
         p{\bf X}(s) =
           {\bf A}(s)~{\bf X}(s) + {\bf B}(s)~W(s), \nonumber \\
         Y(s) = {\bf C}(s)~{\bf X}(s),
   \hspace{3.3em} t \geq 0. \hspace{0.2em}
\end{eqnarray}
where matrix ${\bf A}(s)$ and vectors ${\bf B}(s)$, ${\bf C}(s)$ are
\begin{displaymath} \label{r6}
{\bf A}(s)=               \hspace{-0.2em}
\left[                    \hspace{-0.4em}
 \begin{array}{cc}
    0    &     1   \\
    \,    &    \,  \\
   -\frac{a_0+K}{a_2}s^{2-{\alpha}}&
   -\frac{a_1s^{1+{\beta}-{\alpha}}+T_ds^{1+{\delta}-{\alpha}}}{a_2}
 \end{array}              \hspace{-0.4em}
\right]\hspace{-0.2em},
{\bf B}(s)=               \hspace{-0.2em}
\left[                    \hspace{-0.4em}
 \begin{array}{c}
    0          \\
    \,        \\
    \frac{1}{a_2}s^{2-{\alpha}}
 \end{array}              \hspace{-0.4em}
\right]\hspace{-0.2em},
{\bf C}(s)=               \hspace{-0.2em}
\left[                    \hspace{-0.4em}
 \begin{array}{c}
    K          \\
    \,        \\
    T_d s^{{\delta}-1}
 \end{array}              \hspace{-0.4em}
\right]^T
\end{displaymath}
From equations (\ref{r5}) we can derive the
overall transfer function
\begin{equation} \label{r7}
       \frac{Y(s)}{W(s)} =
       \frac{{\bf C} \, adj(s{\bf I}-{\bf A}) \, {\bf B}}
            {det(s{\bf I} - {\bf A})}
\end{equation}
The characteristic equation of the closed-loop system is determined by
solving the determinant in the denominator of the transfer function,
which is a fractional-order polynomial in $s$
\begin{equation} \label{r8}
       a_2 s^{\alpha} + a_1 s^{\beta} + T_d s^{\delta} + (a_0 + K) = 0
\end{equation}
Solving this
equation we can find the
poles of the closed-loop system with known parameters, or design the
$PD^\delta$ controller parameters for desired control system poles.

Consider the system described by differential equation (\ref{r2})
with system coefficients
$a_2 = 0.8, a_1 = 0.5, a_0 = 1, \alpha = 2.2, \beta = 0.9$.
The task is to determine the controller parameters
$K, T_d$ and $\delta$ for desired control system poles
$s_{1,2}=-1{\pm}6i$ and the steady state error less than 4\%.
Parameter $K=24$ can be computed from the equation
based on the steady state error $e_{ss}(\infty)$ of the closed-loop
control system
$K = ( 100 / e_{ss} - 1 ) a_0$.
If we introduce the desired system poles $s_{1,2}$ to the characteristic
equation (\ref{r8}) we obtain system of two nonlinear equations from
which we can easy derive two equations
$\delta = \arctan(2,9839) \frac{1,8098}{\pi}$  and
${T_d}^2 37^\delta = 645,2174$ \, for the calculation of the
controller parameters $T_d = 6,9407$ and  $\delta = 0,71859$.
In Fig.2 and Fig.3 are depicted the unit - step responses
of the classical numerical solution
\cite{Dorcak2}
and the numerical solution of the state space model (\ref{rrr4}).
The obtained state trajectory represents stable focal point
for the above-mentioned coefficients of the system.

\vspace *{3mm}
\begin{tabular}{cc}
\includegraphics[scale=0.3]{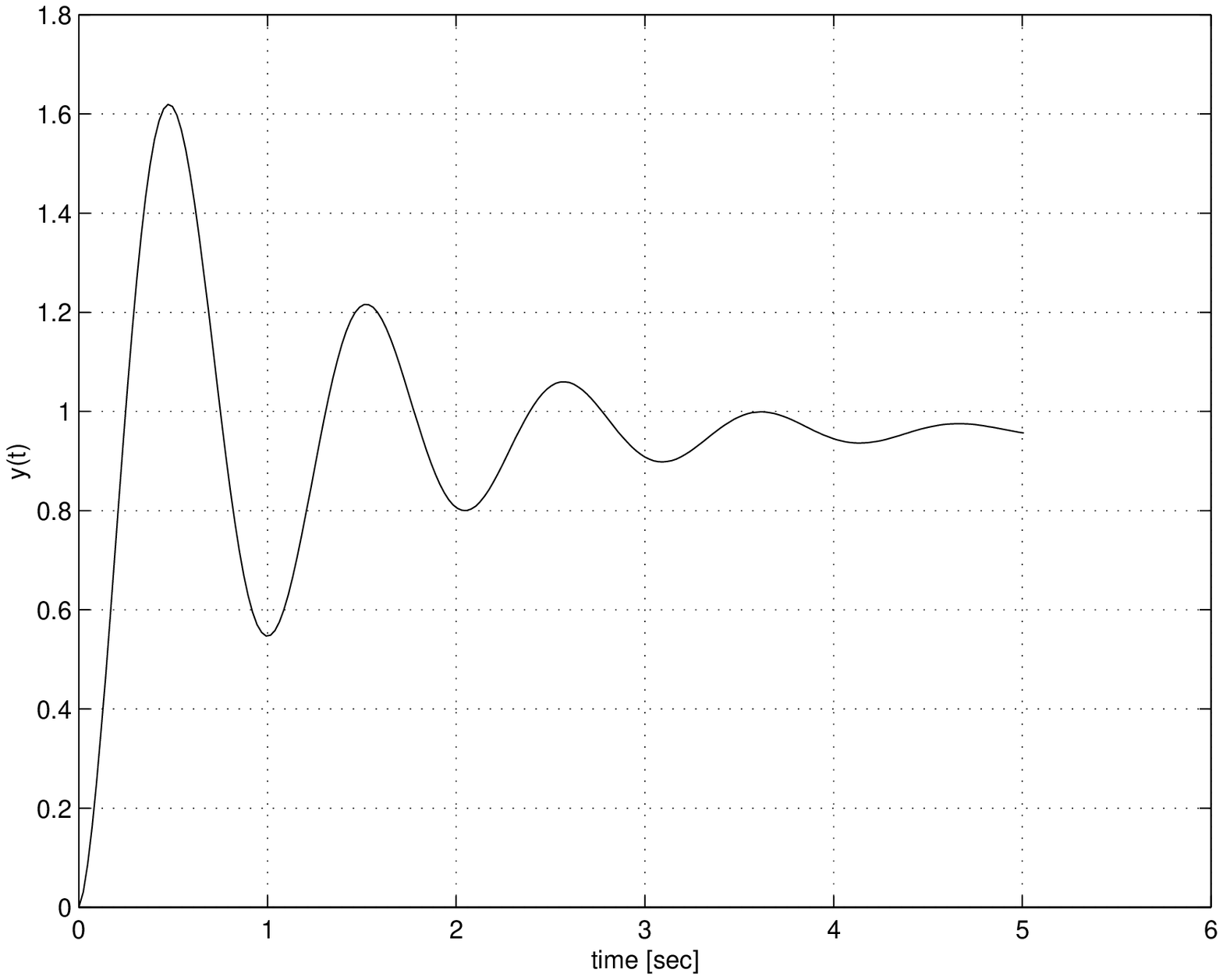} &
\hspace {2mm}
\includegraphics[scale=0.3]{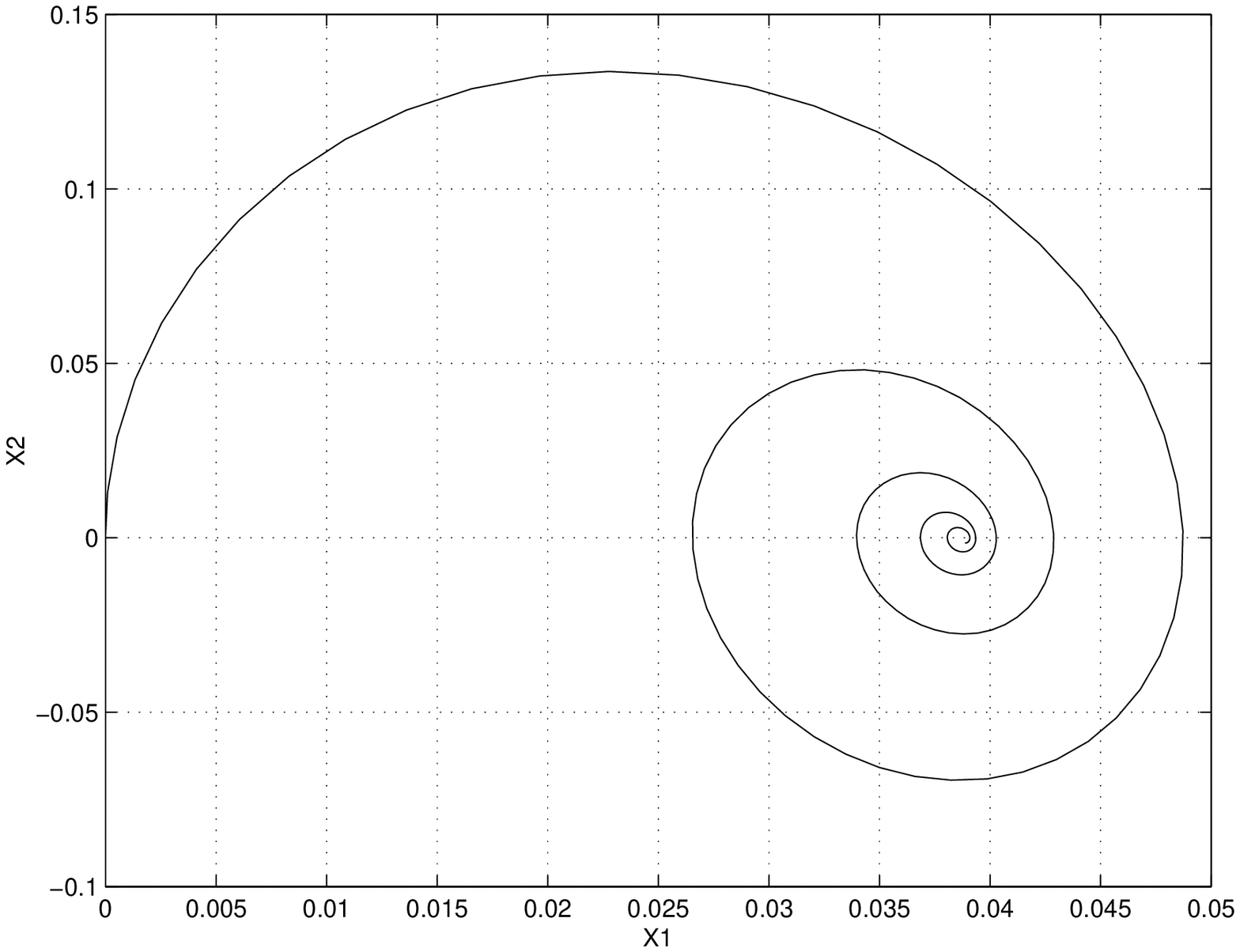} \\ 
\vspace *{2mm}
Figure 2: Unit step response &
Figure 3: State trajectory   \\ 
\end{tabular}

\setcounter{figure}{3}

Assuming the integer-order $PD$ controller we obtain a system of two
linear equations for direct calculation of the controller parameters
$K = 36,0854$ and  $T_d = 4,0141$.

If we increased in case of fractional-order $PD^\delta$ controller
the requirement on steady state error less than 2\%
we would obtain the following controller parameters
$K = 49$, $T_d = -79,74427$ and  $\delta = -0,55194$. That means, we
obtained the controller with weak integrator. But the characteristic
equation of such a system has one additional pole
$s_3 = 1,98$
and it follows from the stability analysis that such closed-loop
system is unstable. We can verify this fact also with frequency methods
for stability investi-gation of the fractional-order system
\cite{Dorcak5, PetrasDorcak1}. In this case we can change the
desired control system poles $s_{1,2}$, or we have to proceed
in a different way, e.g. as follows.

\section*{\normalsize{5. \hspace{1.2mm}Fractional-order~control~system~with
                         $PI^{\lambda}$~controller~in~the~state~space}}

\hspace{5.0mm}
        The differential equation of the closed regulation
system for the transfer function of the controlled system
$G_{s}(s)=1/(a_2s^{\alpha}+a_1s^{\beta}+a_0)$
and the fractional-order $PI^{\lambda}$ controller
$G_{r}(s)=K+T_is^{-\lambda}$
has the form
\vspace{-0mm}
\begin{equation} \label{r10}
     a_{2}\, y^{(\alpha+\lambda)}(t) +
     a_{1}\, y^{(\beta+\lambda)}(t) +
     (a_{0}+K)\, y^{(\lambda)}(t) +
     T_i\, y(t) =
     K\, w^{(\lambda)}(t) +
     T_i\, w(t)
\end{equation}
where  $\alpha, \beta, \lambda$ are generally real
numbers and $a_0, a_1, a_2, K, T_i$ are arbitrary
constants.

\hspace{-7.0mm}
The state space model equivalent to model (\ref{r10}) has
the folloving form
\vspace{-2mm}
\begin{eqnarray} \label{r11}
    \dot{x_1}(t) = - x_2(t) + w ,
     \hspace{9.2cm}
     \nonumber \\
    \dot{x_2}(t) = + x_3(t)~,
     \hspace{9.8cm}
     \nonumber \\
    \dot{x_3}(t) = + \frac{T_i}  {a_2}~x_1^{(3-\alpha-\lambda)}(t)
                   - \frac{a_0+K}{a_2}~x_2^{(2-\alpha)}(t)
                   - \frac{a_1}  {a_2}~x_1^{(1+\beta-\alpha)}(t)
                   + \frac{K}    {a_2}~w_1^{(2-\alpha)}(t)~,
     \hspace{0em} \\
        y  (t) = x_2(t)~,
     \hspace{1.5em} t \geq 0. \hspace{24.8em}
     \nonumber
\vspace{-2mm}
\end{eqnarray}
where $\dot{x_1}(t) = \dot{z}(t)$ is the
actuating error signal and
$\dot{x_2}(t) = \dot{v_1}(t), \dot{x_3}(t) = \dot{v_2}(t)$ are
the state variables of the controlled system (Fig.4).
\begin{figure}[ht]
\centering
\includegraphics[scale=0.5]{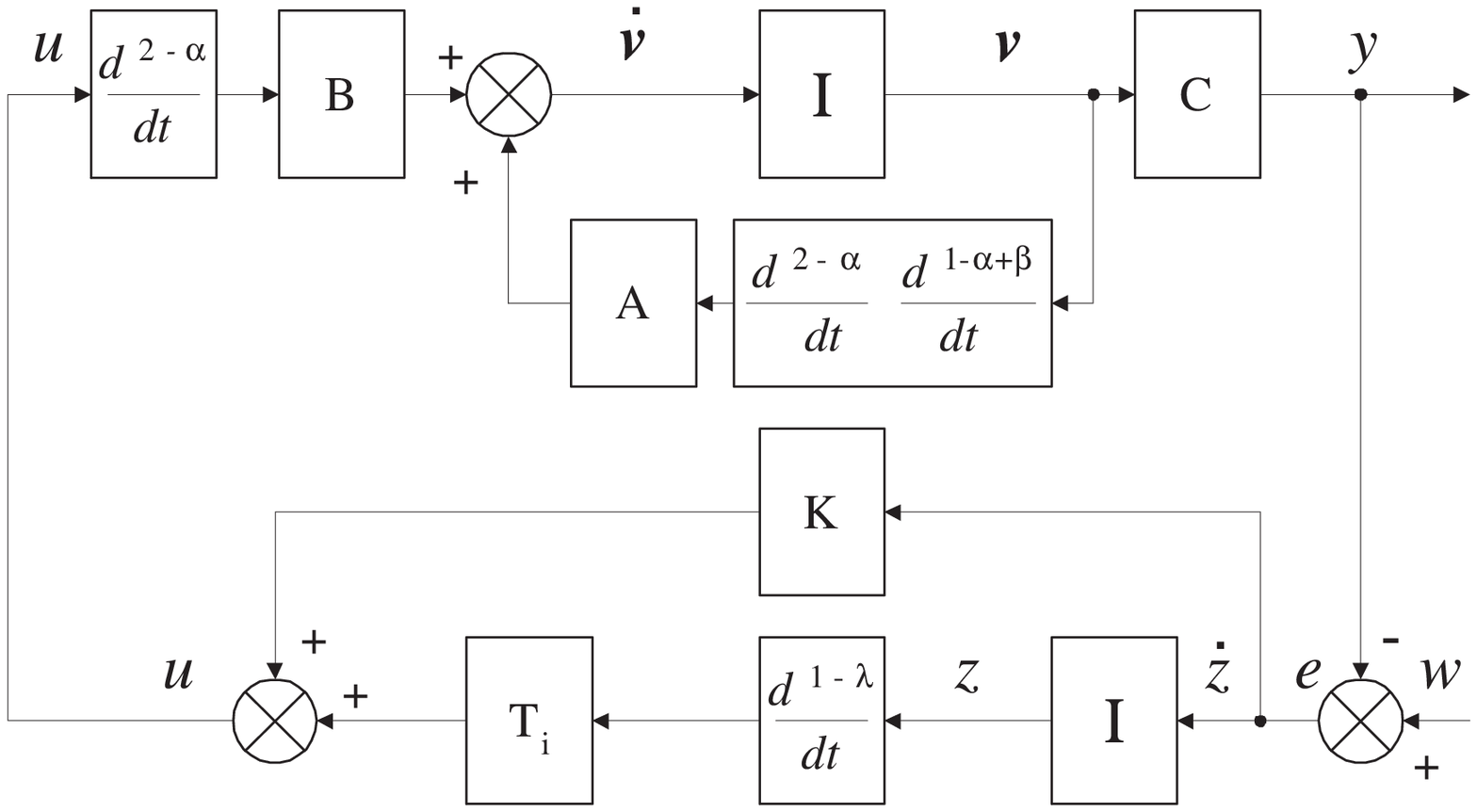}
\caption{Block diagram representing the state and output equations}
         \label{fig:PI}
\end{figure}
As above, we can derive the following characteristic equation
of the closed-loop system
\begin{equation} \label{r12}
       a_2 s^{\alpha+\lambda} + a_1 s^{\beta+\lambda} +
       (a_0 + K) s^{\lambda} + T_i = 0
\end{equation}
Assume the same controlled system as in the previous section.
%
%
If we introduce the desired system poles $s_1$, $s_2$, $s_3$
to the characteristic equation (\ref{r12}) we obtain
a system of three nonlinear equations from which we can
numericaly calculate the controller parameters
$K$, $T_i$ and  $\lambda$.

\vspace{-2.5mm}
\section*{\normalsize{6. \hspace{1.2mm}Conclusion}}
\vspace{-1.5mm}
\hspace{5.0mm}
         We have presented a mathematical description of a
fractional-order control system in the state space and the
state-space controller design based on placing the closed-loop
poles on the complex plane.
In the design of a fractional-order $PD^\delta$ regulator it is
necessary to pay attention to the question of stability of the
control system. In the case of negative values of $\delta$ the
order of the system increases. In addition to the required poles,
the system can develop
new
poles that can render the system unstable.

\vspace{-2.5mm}
\section*{\normalsize{Acknowledgements}}
\vspace{-1.5mm}
\hspace{5.0mm}
This work was partially supported by grant VEGA 1/7098/20 from the
Slovak Grant Agency for Science.

\end{document}